# FOURIER COEFFICIENTS OF HALF-INTEGRAL
# WEIGHT MODULAR FORMS MODULO $\ell$

KEN ONO AND CHRISTOPHER SKINNER

*Dedicated to the memory of S. Chowla.*

## 1. INTRODUCTION

S. Chowla conjectured that every prime $p$ has the property that there are infinitely many imaginary quadratic fields whose class number is not a multiple of $p$. Gauss' genus theory guarantees the existence of infinitely many such fields when $p = 2$, and the work of Davenport and Heilbronn [D-H] suffices for the prime $p = 3$. In addition, the Davenport-Heilbronn result demonstrates that a positive proportion of such fields have this property. Using an elementary argument based on the Kronecker recurrence relations, Hartung [Ha] proved that for any odd prime $p$ there are indeed infinitely many such fields whose class numbers are not multiples of $p$. His argument has been employed in other similar studies [Ho1,Ho2, Ho-On].

It is easy to capture the flavor of the Kronecker relations: if $r(n)$ denotes the number of representations of a positive integer $n$ as a sum of three squares, then

$$(1) \qquad \sum_{k=-\infty}^{\infty} r(n - k^2) = 8 \sum_{\substack{d \mid n \\ d \not\equiv 0 \pmod 4}} d.$$

By Gauss' theorem, (1) is a recurrence relation for class numbers since

$$r(n) = \begin{cases} 12h(-4n) & \text{if } n \equiv 1, 2, 5, 6 \pmod 8, \\ 24h(-n) & \text{if } n \equiv 3 \pmod 8 \end{cases}$$

for square-free $n > 3$. If $\Theta(z) := 1 + 2q + 2q^4 + 2q^9 + 2q^{16} + \dots$ ($q := e^{2\pi i z}$ throughout), then (1) is simply the fact that the product of $\Theta^3(z)$, the 'generating function' for $r(n)$, and $\Theta(z)$ is a normalized Eisenstein series of weight 2 on $\Gamma_0(4)$.

In light of the work of Shimura [Sh1] and Waldspurger [Wal], it is now desirable to study the indivisibility properties of the coefficients of half-integral weight modular forms. Such results shed light on the structure of Tate-Shafarevich groups of elliptic curves, Bloch-Kato III's of newforms, class numbers of quadratic fields, and generalized Bernoulli numbers. In earlier work, Jochnowitz [Jo] obtained indivisibility results for such coefficients by developing a theory of half-integral weight modular forms mod $\ell$ analogous to the integral

1991 *Mathematics Subject Classification.* Primary 11F33; Secondary 11G40.

The first author is supported by NSF grants DMS-9508976 and DMS-9304580, and the second author is supported by a Sloan Foundation Doctoral Dissertation Fellowship.

Typeset by $\mathcal{A}\mathcal{M}\mathcal{S}$-TEX

weight theory developed by Serre and Swinnerton-Dyer. We show that such results can be obtained using Kronecker's idea of multiplying by theta functions. This approach has the pleasant property that it uncovers congruences like (18) which are mod $\ell$ analogs of the Kronecker relations.

Let $f(z) = \sum_{n=1}^{\infty} a(n)q^n \in S_{2k}(N, \chi_0)$ be a newform with trivial Nebentypus $\chi_0$, and let $L(f, s) := \sum_{n=1}^{\infty} a(n)n^{-s}$ be its $L$-function. Then $\Lambda(f, s) := (2\pi)^{-s}\Gamma(s)N^{s/2}L(f, s)$ satisfies the functional equation

$$\Lambda(f, s) = \epsilon \cdot \Lambda(f, 2k - s),$$

where $\epsilon = \pm 1$ is the sign of the functional equation. We shall be interested in the quadratic twists of this function. If $D \neq 0$, then let $\chi_D$ denote the Kronecker character for the field $\mathbb{Q}(\sqrt{D})$, and if $D = 0$, then let $\chi_0$ denote the trivial character. For notational convenience, if $D$ is a fundamental discriminant of a quadratic number field, then define $D_0$ by

(2)
$$D_0 := \begin{cases} |D| & \text{if } D \text{ is odd,} \\ |D|/4 & \text{if } D \text{ is even.} \end{cases}$$

The $D$−quadratic twist of $f$, denoted $f_D$, is the newform corresponding to the twist of $f$ by the character $\chi_D$. If $D = 0$ or if $(D, N) = 1$, then $f_D(z) := \sum_{n=1}^{\infty} \chi_D(n)a(n)q^n$. Its $L$-function satisfies the functional equation

$$\Lambda(f_D, s) = \epsilon \cdot \chi_D(-N)\Lambda(f_D, 2k - s).$$

Waldspurger proved [Th. 1, Wal] that for suitable $f$, the central critical values $L(f_D, k)$ for suitable fundamental discriminants $D$ where $(-1)^k D > 0$ are typically given by

$$L(f_D, k) = c(D_0)^2 \cdot A_D$$

where $A_D \neq 0$ and $g(z) = \sum_{n=1}^{\infty} c(n)q^n$ is a half integral weight cusp form whose Shimura lift is $f(z)$. In view of the number theoretic significance of central critical values of $L$−functions, it is natural to study the indivisibility properties of the coefficients $c(D_0)$. For example consider Kolyvagin's Conjecture [Conjecture F, Ko2]. If $E/\mathbb{Q}$ is an elliptic curve, then let $\text{III}(E)$ denote its Tate-Shafarevich group, $L(E, s)$ its Hasse-Weil $L$−function, and $\Omega_E$ its period as defined in §2. Furthermore let $E_D$ denote the $D$−quadratic twist of $E$.

**Kolyvagin's Conjecture.** *Let $E/\mathbb{Q}$ be an elliptic curve, $\ell$ a prime, and $S(E, D)$ the conjectured order of $\text{III}(E_D)$.*

(i) *If $L'(E, 1) \neq 0$, then there is a fundamental discriminant $D$ for which*

$$L(E_D, 1) \neq 0 \quad \text{and} \quad S(E, D) \not\equiv 0 \pmod{\ell}.$$

(ii) *If $L(E, 1) \neq 0$, then there is a fundamental discriminant $D$ for which*

$$L'(E_D, 1) \neq 0 \quad \text{and} \quad S(E, D) \not\equiv 0 \pmod{\ell}.$$

For further implications of this conjecture see [B-D].

For exceptional $g(z) = \sum_{n=1}^{\infty} c(n)q^n$, explicit relations like (1) can be obtained where divisor functions are replaced by values of Hecke Grössencharacters. For instance in [O] such relations are exhibited for the weight $\frac{3}{2}$ cusp forms pertaining to the Tate-Shafarevich groups $\mathrm{III}(N)$ of the elliptic curves over $\mathbb{Q}$ of the form $E(N): y^2 = x^3 - N^2 x$. Assuming the Birch and Swinnerton-Dyer Conjecture, if $2N-1$ is a positive square-free integer for which $E(2N-1)$ has rank 0, then

$$\mathfrak{T}(2N-1,1)\sqrt{|\mathrm{III}(2N-1)|} + \sum_{(m,k)\in\mathfrak{S}(N)} \mathfrak{T}\left(\frac{2N-k^2}{m^2},m\right)\sqrt{\left|\mathrm{III}\left(\frac{2N-k^2}{m^2}\right)\right|}$$

$$= \sum_{(x,y)\in\mathfrak{F}(N)} (-1)^y(2y+1) + 2\sum_{(x,y)\in\mathfrak{F}(N/2)} (-1)^{x+y}(2y+1).$$

Here $\mathfrak{S}(N)$ and $\mathfrak{F}(N)$ are explicit finite sets, and $\mathfrak{T}(t,m)$ is an explicit function involving Jacobsthal sums. Using a theorem of Rubin [Ru], it is then easy to deduce that if $\ell$ is a prime, then there are infinitely many odd square-free integers $N$ for which $E(N)$ has rank 0 and $|\mathrm{III}(N)| \not\equiv 0 \pmod{\ell}$.

However, such relations are very rare. Fortunately, one can use $\ell-$adic Galois representations to prove a general theorem yielding indivisibility results as strong as those obtained by such explicit relations. In §2 we prove such a result and obtain:

**Corollary 1.** *Let $f(z) \in S_{2k}(N, \chi_0)$ be a newform, and let $\Omega$ be a period for $f$. Then for all but finitely many primes $\ell$ there are infinitely many fundamental discriminants $D$ for which*

$$(-1)^k D > 0 \quad \text{and} \quad \left|\frac{L(f_D, k)D_0^{k-\frac{1}{2}}}{\Omega}\right|_\ell = 1.$$

We obtain the following partial resolution to part $(i)$ of Kolyvagin's conjecture.

**Corollary 2.** *Let $E/\mathbb{Q}$ be a modular elliptic curve. If $S(E, D)$ denotes the order of $\mathrm{III}(E_D)$ as predicted by the Birch and Swinnerton-Dyer Conjecture, then for all but finitely many primes $\ell$ there are infinitely many fundamental discriminants $D < 0$ for which*

$$L(E_D, 1) \neq 0 \quad \text{and} \quad S(E, D) \not\equiv 0 \pmod{\ell}.$$

**Corollary 3.** *If $E/\mathbb{Q}$ is an elliptic curve with complex multiplication, then for all but finitely many primes $\ell$ there are infinitely many fundamental discriminants $D < 0$ for which $L(E_D, 1) \neq 0$ and*

$$|\mathrm{III}(E_D)| \not\equiv 0 \pmod{\ell}.$$

We also obtain a 'local' version of Corollary 1.

**Corollary 4.** *Let $f(z) \in S_{2k}(N, \chi_0)$ be a newform, and let $\Omega$ be a period for $f$. Let $P := \{p_1, p_2, \ldots p_s\}$ be a finite set of primes, and $\epsilon := (\epsilon_1, \epsilon_2, \ldots \epsilon_s)$ where $\epsilon_i \in \{\pm 1\}$. For all but finitely many primes $\ell$ there are infinitely many fundamental discriminants $D$ for which $(-1)^k D > 0$, $\chi_D(p_i) = \epsilon_i$ for every $1 \leq i \leq s$, and*

$$\left|\frac{L(f_D, k)D_0^{k-\frac{1}{2}}}{\Omega}\right|_\ell = 1.$$

**Corollary 5.** *Let $E/\mathbb{Q}$ be a modular elliptic curve with conductor $N$ for which $L(E, s)$ has a simple zero at $s = 1$. Then $E$ has rank 1, and for all primes $\ell$ outside a finite set which is effectively determinable (see Remark 4)*

$$ord_\ell(\text{Ш}(E)) \leq ord_\ell(S(E)),$$

*where $S(E)$ denotes the order of $\text{Ш}(E)$ as predicted by the Birch and Swinnerton-Dyer Conjecture.*

## 2. Results

If $k$ is a positive integer, then let $S_k(N)$ denote the space of cusp forms of weight $k$ on $\Gamma_1(N)$, and let $S_k^{cm}(N)$ denote the subspace of $S_k(N)$ spanned by those forms having complex multiplication (cf. [R2]). If $\chi$ is a Dirichlet character modulo $N$, then let $S_k(N, \chi)$ denote the subspace of $S_k(N)$ consisting of those forms having Nebentypus character $\chi$. By the theory of newforms, every $f(z) \in S_k(N)$ can be uniquely expressed as a linear combination

$$f(z) = \sum_{i=1}^{r} \alpha_i A_i(z) + \sum_{j=1}^{s} \beta_j B_j(\delta_j z),$$

where $A_i(z)$ and $B_j(z)$ are newforms of weight $k$ and level a divisor of $N$, and where each $\delta_j$ is a non-trivial divisor of $N$. Let

$$f^{new}(z) := \sum_{i=1}^{r} \alpha_i A_i(z) \quad \text{and} \quad f^{old}(z) := \sum_{j=1}^{s} \beta_j B_j(\delta_j z)$$

be, respectively, the *new part of $f$* and the *old part of $f$*.

For a non-negative integer $k$, let $M_{k+\frac{1}{2}}(N)$ (resp. $S_{k+\frac{1}{2}}(N)$) denote the space of modular forms (resp. cusp forms) of half-integral weight $k + \frac{1}{2}$ on $\Gamma_1(4N)$. If $i = 0$ or $1$, $0 \leq r < t$, and $a \geq 1$, then let $\theta_{a,i,r,t}(z)$ denote the Shimura theta function

$$\theta_{a,i,r,t}(z) := \sum_{n \equiv r \pmod{t}} n^i q^{an^2}.$$

If $\Theta(N)$ is the space of such functions of level $4N$, then the Serre-Stark Theorem [S-S] implies

$$\Theta(N) = M_{\frac{1}{2}}(N) \cup \left\{ \text{subspace of } M_{\frac{3}{2}}(N) \text{ spanned by those } \theta_{a,1,r,t}(z) \text{ on } \Gamma_1(4N) \right\}.$$

If $g(z) \in M_{k+\frac{1}{2}}(N)$ and $h(z) \in \Theta(N')$, then $G_h(z) := g(z) \cdot h(z)$ is a modular form on $\Gamma_1(4NN')$ of integral weight $k + 1$ or $k + 2$.

**Definition 1.** *A modular form $g(z) \in M_{k+\frac{1}{2}}(N)$ is **good** if there exists some $N'$ and some $h(z) \in \Theta(N')$ for which*

(G1)  $G_h(z)$ *is a cusp form.*

(G2)  $0 \neq G_h^{new}(z) \notin S_{k+1}^{cm}(4NN') \cup S_{k+2}^{cm}(4NN')$.

Let $\overline{\mathbb{Q}}$ be an algebraic closure of $\mathbb{Q}$, and for each rational prime $\ell$, let $\overline{\mathbb{Q}}_\ell$ be an algebraic closure of $\mathbb{Q}_\ell$. Fix an embedding of $\overline{\mathbb{Q}}$ into $\overline{\mathbb{Q}}_\ell$. This fixes a choice of decomposition group $D_\ell$. In particular, if $K$ is any finite extension of $\mathbb{Q}$, and if $\mathcal{O}_K$ is the ring of integers of $K$, then for each $\ell$ this fixes a choice of a prime ideal $\mathfrak{p}_{\ell,K}$ of $\mathcal{O}_K$ dividing $\ell$. Let $\mathbb{F}_{\ell,K}$ be the residue field of $\mathfrak{p}_{\ell,K}$, and let $|\cdot|_\ell$ be an extension to $\overline{\mathbb{Q}}_\ell$ of the usual $\ell$-adic absolute value on $\mathbb{Q}_\ell$.

If $f(z) = \sum_{n=1}^\infty a(n)q^n \in S_k(N,\chi)$ is a newform, then the $a(n)$'s are algebraic integers and generate a finite extension of $\mathbb{Q}$, say $K_f$. If $K$ is any finite extension of $\mathbb{Q}$ containing $K_f$, and if $\ell$ is any prime, then by the work of Deligne and Serre [D], [D-S] there is a continuous, semisimple representation

$$\rho_{f,\ell} : Gal(\overline{\mathbb{Q}}/\mathbb{Q}) \to GL_2(\mathbb{F}_{\ell,K})$$

for which

(R1)  $\rho_{f,\ell}$ is unramified at all primes $p \nmid N\ell$.

(R2)  trace $\rho_{f,\ell}(\text{frob}_p) = a(p) \bmod \mathfrak{p}_{\ell,K}$ for all primes $p \nmid N\ell$.

(R3)  det $\rho_{f,\ell}(\text{frob}_p) = \chi(p)p^{k-1} \bmod \mathfrak{p}_{\ell,K}$ for all primes $p \nmid N\ell$.

(R4)  det $\rho_{f,\ell}(c) = -1$  for any complex conjugation $c$.

Here $\text{frob}_p$ denotes any Frobenius element for the prime $p$. These representations capture many properties of the reductions of the $a(n)$'s modulo $\mathfrak{p}_{\ell,K}$. For example, if $\Delta(z) = \sum_{n=1}^\infty \tau(n)q^n$ is Ramanujan's $\Delta$-function, then the representation $\rho_{\Delta,691}$ is

$$\rho_{\Delta,691} \sim \begin{pmatrix} 1 & \\ & \omega^{11} \end{pmatrix}.$$

Here $\omega$ denotes the cyclotomic character giving the action of $Gal(\overline{\mathbb{Q}}/\mathbb{Q})$ on the 691th roots of unity, so $\omega(\text{frob}_p) = p$ for all $p \neq 691$. This 'explains' the classical Ramanujan congruence

$$\tau(n) \equiv \sigma_{11}(n) \pmod{691}$$

by $(R2)$ and the multiplicativity of $\tau(n)$. For most cusp forms and most primes, the corresponding representations do not have such a nice description (see $(L6)$ below), a fact which is exploited in the proof of the following theorem.

**Theorem 1.** *Let $g(z) = \sum_{n=0}^\infty c(n)q^n \in M_{k+\frac{1}{2}}(N)$ be an eigenform whose coefficients are algebraic integers. If $g(z)$ is good, then for all but finitely many primes $\ell$ there exist infinitely many square-free integers $m$ for which*

$$|c(m)|_\ell = 1.$$

*Proof.* Since $g(z)$ is good, there exists an $h(z) \in \Theta(N')$ for which $G_h(z) = g(z)h(z)$ satisfies $(G1)$ and $(G2)$. If $G_h(z) = \sum_{n=1}^\infty a(n)q^n$, and if $h(z) = \sum_{n=0}^\infty b(n)q^{an^2}$, then

$$(3) \qquad\qquad a(n) = \sum_{\substack{x,y \\ x+ay^2=n}} c(x)b(y).$$

Since the $c(x)$'s and $b(y)$'s are algebraic integers, the same is true of the $a(n)$'s.

Let $k'$ be the weight of $G_h(z)$ (so $k'$ is either $k+1$ or $k+2$), and let $\mathcal{N}_0$ be the set of newforms of weight $k'$ and level dividing $4NN'$. The new part of $G_h(z)$ can be uniquely expressed as a linear combination

$$(4) \qquad G_h^{new}(z) = \sum_{f(z) \in \mathcal{N}_0} \alpha_f f(z).$$

Since the Fourier coefficients of $G_h(z)$ are algebraic, the same is true for those of $G_h^{new}(z)$ and $G_h^{old}(z)$, and each $\alpha_f$.

A critical feature of the proof is the 'proper' bookkeeping of newforms which are related via twisting and Galois conjugation. We now fix our notation. Let $X$ be the set of Dirichlet characters of conductor dividing $8NN'$, and let

$$\mathcal{N} = \{f_\chi(z) \ : f \in \mathcal{N}_0, \ \chi \in X\},$$

where $f_\chi$ is the newform corresponding to the twist of $f$ by $\chi$. Each newform in $\mathcal{N}$ has level divisible only by primes dividing $4NN'$ and the conductor of its Nebentypus character is a divisor of $4NN'$. Let $K$ be a finite, Galois extension of $\mathbb{Q}$ containing the Fourier coefficients of each $f \in \mathcal{N}$, the $\alpha_f$'s, the Fourier coefficients of $G_h^{old}(z)$, and a primitive $8NN'$th root of unity.

If $f \in \mathcal{N}$, then let $G_f \subseteq Gal(\overline{\mathbb{Q}}/\mathbb{Q})$ be the subgroup stabilizing the set

$$\mathcal{N}_f := \{f_\chi(z) \ : \ \chi \in X\}.$$

For each prime $\ell$, let $D_{f,\ell} := G_f \cap D_\ell$, and let $\mathbb{F}_{f,\ell} := \mathbb{F}_{\ell,K}^{D_{f,\ell}}$.

Fix a prime $\ell$ for which:

$(L1)$  $K$ is unramified at $\ell$.

$(L2)$  $\ell \nmid 4NN'$ and $\ell > 2(k+2)$.

$(L3)$  If $\alpha_f \neq 0$, then $|\alpha_f|_\ell = 1$.

$(L4)$  The characters $\chi \in X$, viewed as taking values in $\mathbb{F}_{\ell,K}$, are distinct.

$(L5)$  The representations $\rho_{f,\ell}$  $(f \in \mathcal{N})$ are pairwise non-isomorphic.

$(L6)$  If $f \in \mathcal{N}$ does not have complex multiplication, then the image of $\rho_{f,\ell}$ contains a normal subgroup $H_f$ conjugate to $SL_2(\mathbb{F}_{f,\ell})$ and for which $\mathrm{Im}\rho_{f,\ell}/H_f$ is abelian.

It is clear that $(L1)-(L4)$ hold for all sufficiently large primes $\ell$, and it follows from $(R2)$ and 'multiplicity one' for newforms that $(L5)$ is also true for all large primes. Ribet [R3], following the ideas of Serre [Se], has shown that $(L6)$ holds for all large primes.

For each $f \in \mathcal{N}_0$, let $S_f \subset D_\ell$ be a set of representatives for the classes $D_\ell/D_{f,\ell}$, and put $\mathcal{N}_{f,\ell} = \cup_{\sigma \in S_f} \mathcal{N}_{f^\sigma}$, a disjoint union. Let $\mathcal{N}' = \{f_1(z), f_2(z), ..., f_u(z)\}$ be a subset of $\mathcal{N}_0$ such that $\mathcal{N}$ is a disjoint union of the $\mathcal{N}_{f_i,\ell}$'s. We assume that $f_1(z), ..., f_v(z)$ do not have complex multiplication and that $f_{v+1}(z), ..., f_u(z)$ do. Since $G_h(z)$ satisfies $(G2)$, $f_1(z)$ can be chosen so that the coefficient $\alpha_{f_1}$ of $f_1(z)$ in (4) is non-zero. Write $f_i(z) = \sum_{n=1}^{\infty} a_i(n)q^n$, write $\rho_i$ for $\rho_{f_i,\ell}$, $S_i$ for $S_{f_i}$, $H_i$ for $H_{f_i}$, and $\mathbb{F}_i$ for $\mathbb{F}_{f_i,\ell}$. In this way we are able to conveniently keep track of those distinct newforms which are closely related.

**Lemma 1.** *With the above notation,*

(i) *The image of $\rho_1 \times \cdots \times \rho_v$ contains a normal subgroup conjugate to $SL_2(\mathbb{F}_1) \times \cdots \times SL_2(\mathbb{F}_v)$.*

(ii) *If $f_i(z)$ does not have complex multiplication, then for each positive integer $d$ and each $w \in \mathbb{F}_i$, a positive density of primes $p \equiv 1 \pmod{d}$ satisfies*

$$a_i(p) \equiv w \pmod{\mathfrak{p}_{\ell,K}}.$$

(iii) *If $f_i(z)$ does not have complex multiplication, then for each pair of co-prime positive integers $r, d$, a positive density of primes $p \equiv r \pmod{d}$ satisfies*

$$|a_i(p)|_\ell = 1.$$

*Proof of Lemma 1.* Part $(i)$ is surely well-known, and in any event can be proved by a simple modification of the arguments in [R1], [R3], and [Se]. Without loss of generality, we can assume that $H_i = SL_2(\mathbb{F}_i)$. Let

$$H := \bigcap_{i=1}^{v} \rho_i^{-1}(H_i).$$

The quotient $H_i/\rho_i(H)$ is abelian by $(L6)$, and since $SL_2(\mathbb{F}_i)$ has no non-trivial abelian quotients, $\rho_i(H) = H_i$. It suffices to show that $(\rho_1 \times \cdots \times \rho_v)(H)$ contains $SL_2(\mathbb{F}_1) \times \cdots \times SL_2(\mathbb{F}_v)$. By [Lemma 3.3, R1], it suffices to show that $(\rho_i \times \rho_j)(H)$ contains $SL_2(\mathbb{F}_i) \times SL_2(\mathbb{F}_j)$ if $i \neq j$. Arguing as in the proofs of [Lemme 8, §6.2, Se] and [Theorem 3.8, R1], one easily sees that if this is not true, then there is a $\sigma \in D_\ell$ such that $\rho_i^\sigma \cong \rho_j \otimes \phi$ for some Dirichlet character $\phi$ unramified at primes not dividing $4NN'\ell$. Arguing as in the proof of [Theorem 6.1, R1], one sees that since $\ell > 2(k+2)$ (see $(L2)$) $\phi$ is unramified at $\ell$. Ribet's arguments employ various results of Swinnerton-Dyer about level 1 modular forms mod $\ell$, which have been generalized to higher level by Katz [Ka] (see also [§4, Gr]). Now, comparing determinants shows that the conductor of $\phi^2$ divides $4NN'$, so the conductor of $\phi$ divides $8NN'$. Therefore, $f^\sigma(z)$ and $f_\phi(z)$ are in $\mathcal{N}$ and $\rho_{f^\sigma,\ell} \cong \rho_{f_\phi,\ell}$, contradicting $(L5)$. This proves part $(i)$.

Parts $(ii)$ and $(iii)$ are simple consequences of property $(L6)$. Let $\epsilon : Gal(\overline{\mathbb{Q}}/\mathbb{Q}) \to (\mathbb{Z}/d)^\times$ be the character defined by $\epsilon(\sigma) = s$ if $\zeta_d^\sigma = \zeta_d^s$. The condition $p \equiv r \pmod{d}$ is equivalent to $\epsilon(\mathrm{frob}_p) = r$. Since $SL_2(\mathbb{F}_i)$ has no non-trivial abelian quotient, the image of $\rho_i \times \epsilon$ contains $SL_2(\mathbb{F}_i) \times 1$. Let $g = \begin{pmatrix} 0 & 1 \\ -1 & w \end{pmatrix} \in SL_2(\mathbb{F}_i)$. By the Chebotarev Density Theorem, a positive proportion of primes $p$ satisfy $(\rho_i \times \epsilon)(\mathrm{frob}_p) = (g, 1)$. For such a $p$,

$$a_i(p) \equiv \mathrm{trace} \begin{pmatrix} 0 & 1 \\ -1 & w \end{pmatrix} \equiv w \pmod{\mathfrak{p}_{\ell,K}} \quad \text{and} \quad p \equiv \epsilon(\mathrm{frob}_p) \equiv 1 \pmod{d}.$$

This proves part $(ii)$. Now choose $\sigma \in Gal(\overline{\mathbb{Q}}/\mathbb{Q})$ such that $\epsilon(\sigma) = r$. Clearly, there exists some $g \in SL_2(\mathbb{F}_i)$ for which trace $(\rho_i(\sigma) \cdot g) \neq 0$. Again by the Chebotarev Density Theorem, a positive proportion of the primes $p$ satisfy $(\rho_i \times \epsilon)(\mathrm{frob}_p) = (\rho_i(\sigma), r)(g, 1) = (\rho_i(\sigma) \cdot g, r)$. This proves part $(iii)$.

Q.E.D. Lemma 1

Returning to the proof of the theorem, rewrite (4) as

$$(5) \qquad G_h^{new}(z) = \sum_{i=1}^{u} \sum_{\sigma \in S_i} \sum_{\chi \in X} \alpha(\sigma, i, \chi) f_{i,\chi}^{\sigma}(z),$$

where $\alpha(\sigma, i, \chi) = \alpha_{f_{i,\chi}^{\sigma}}$, and define $a(n)^{new}$ by $G_h^{new}(z) = \sum_{i=1}^{\infty} a(n)^{new} q^n$. If $n$ is a positive integer relatively prime to $4NN'$, then by (5)

$$(6) \qquad a(n)^{new} = \sum_{i=1}^{u} \sum_{\sigma \in S_i} \left( \sum_{\chi \in X} \alpha(\sigma, i, \chi) \chi^{\sigma}(n) \right) a_i(n)^{\sigma}.$$

If $\sigma_1$ is the element of $S_1$ for which $f_1^{\sigma_1} = f_1$, then $\alpha(\sigma_1, 1, \chi_0) = \alpha_{f_1} \neq 0$. By $(L3)$, $|\alpha(\sigma_1, 1, \chi_0)|_{\ell} = 1$, and from $(L4)$, one sees easily that there is an integer $r$ relatively prime to $4NN'$ for which $|\sum_{\chi \in X} \chi^{\sigma_1}(r)|_{\ell} = 1$. Fix such an $r$, and put

$$\lambda(\sigma, i) = \sum_{\chi \in X} \alpha(\sigma, i, \chi) \chi^{\sigma}(r), \quad (1 \leq i \leq u; \ \sigma \in S_i).$$

By (6),

$$a(n)^{new} = \sum_{i=1}^{u} \sum_{\sigma \in S_i} \lambda(\sigma, i) a_i(n)^{\sigma} \quad \text{if} \quad n \equiv r \pmod{8NN'}.$$

Suppose that there are only finitely many square-free intgers $m$ for which $|c(m)|_{\ell} = 1$. Let these be $m_1, ..., m_t$. Let $N''$ be the product of the odd prime factors of $NN'm_1 \cdots m_t$. Let $r_1$ and $d_1$ be relatively prime positive integers such that if $p$ is any prime congruent to $r_1$ modulo $d_1$, then $p$ does not split in any imaginary quadratic subfield having discriminant a divisor of $4N''$. Choose a prime $p_1$ for which $p_1 \nmid 4N'$, $p_1 \equiv r_1 \pmod{d_1}$, and $|a_1(p_1)|_{\ell} = 1$. This is possible by part $(iii)$ of Lemma 1. If $v + 1 \leq i \leq u$, then $a_i(p_1) = 0$ since $f_i(z)$ has complex multiplication. Choose a prime $p_2$ not dividing $4N''$ and for which $|a_1(p_2)|_{\ell} = 1$ and $|a_i(p_2)|_{\ell} < 1$ $(i = 2, ..., v)$. This is possible by part $(i)$ of Lemma 1. For example, let $g \in \text{Im}(\rho_1 \times \cdots \times \rho_v)$ be conjugate to

$$\begin{pmatrix} 1 & 0 \\ 0 & 1 \end{pmatrix} \times \begin{pmatrix} 0 & -1 \\ 1 & 0 \end{pmatrix} \times \cdots \times \begin{pmatrix} 0 & -1 \\ 1 & 0 \end{pmatrix}.$$

By the Chebotarev Density Theorem, a positive proportion of primes $p \nmid 4N''$ satisfy $(\rho_1 \times \cdots \times \rho_v)(\text{frob}_p) = g$. Any such prime clearly has the desired properties. Next, choose a prime $p_3$ for which $p_3 \equiv r(p_1 p_2)^{-1} \pmod{8NN'}$ and $|a_1(p_3)|_{\ell} = 1$. This is possible by part $(iii)$ of Lemma 1.

For $w \in \mathcal{O}_K$, write $\overline{w}$ for its image in $\mathbb{F}_{\ell,K}$. Since $|\lambda(\sigma_1, 1) a_1(p_1 p_2 p_3)^{\sigma_1}|_{\ell} = 1$, and since $|\lambda(\sigma, 1) a_1(p_1 p_2 p_3)^{\sigma}|_{\ell} \leq 1$ for all $\sigma \in S_1$, one sees easily that there is some $w \in \mathcal{O}_K$ for which

$$(7) \qquad \overline{w} \in \mathbb{F}_1^{\times} \quad \text{and} \quad |\sum_{\sigma \in S_1} \lambda(\sigma, 1) a_1(p_1 p_2 p_3)^{\sigma} w^{\sigma}|_{\ell} = 1.$$

Now choose a prime $p_4$ for which $p_4 \equiv 1 \pmod{8NN'}$ and $a_1(p_4) \equiv w \pmod{\mathfrak{p}_{\ell,K}}$. This is possible by part $(ii)$ of Lemma 1. Let $m = p_1 p_2 p_3 p_4$. By the multiplicitivity of the Fourier coefficients of a newform, it follows from the choices of $p_1$ and $p_2$ that

$$|a_i(n)|_{\ell} < 1, \quad (i = 2, ..., u).$$

Therefore, by (7)

$$\left|\sum_{i=2}^{u}\sum_{\sigma\in S_i}\lambda(\sigma,i)a_i(m)^\sigma\right|_\ell < 1 \quad \text{and} \quad \left|\sum_{\sigma\in S_1}\lambda(\sigma,1)a_1(m)^\sigma\right|_\ell = 1,$$

so

$$|a(m)^{new}|_\ell = 1.$$

Since $m$ is prime to $4NN'$, it follows that

(8) $$|a(m)|_\ell = |a(m)^{new}|_\ell = 1.$$

By (3) and the fact that $g(z)$ is an eigenform, $m$ must be of the form

(9) $$m = m_j x^2 + ay^2$$

for some $j \in \{1, ..., t\}$ and some integers $x$ and $y$. Since $m_j X^2 + aY^2$ is a positive definite quadratic form with discriminant, say $d_j$, a divisor of $4N''$, a necessary condition for a solution of (9) is that $\left(\frac{d_j}{p}\right) = 1$ for every prime divisor of $p$ of $m$. But $m = p_1 p_2 p_3 p_4$ was chosen so that $\left(\frac{d_j}{p_1}\right) = -1$ for each $d_j$. This contradiction proves the theorem.

Q.E.D. Theorem 1

This theorem has many arithmetic consequences. Suppose $f(z) \in S_{2k}(N, \chi_0)$ is a newform. If $D$ is a fundamental discriminant, then let $D_0$ be as defined in (2). A non-zero complex number $\Omega \in \mathbb{C}^\times$ is a *period for $f(z)$* if

(10) $$\frac{L(f_D, k)D_0^{k-\frac{1}{2}}}{\Omega}, \quad (-1)^k D > 0,$$

is always an algebraic number (this is a slight abuse of standard terminology). A period $\Omega$ is *nice* if the quantity (10) is always an algebraic integer. That nice periods exist is essentially a result of Shimura. Any period is a multiple of a nice period by some algebraic number. As a consequence of Theorem 1 we obtain Corollary 1.

*Proof of Corollary 1.* Since any period is the multiple of a nice one by an algebraic number, it suffices to prove the corollary for $\Omega$ a nice period. There is a twist $f_\chi$, $\chi$ possibly trivial, satisfying Hypotheses H1 and H2 of [Wal]. By [Th. 1, Wal] there is an integer $N'$ and an eigenform $g(z) = \sum_{n=1}^\infty c(n)q^n \in S_{k+\frac{1}{2}}(N')$ such that for each fundamental discriminant $D$ for which $(-1)^k D > 0$

(11) $$c(D_0)^2 = \begin{cases} \pm \frac{L(f_D, k)D_0^{k-\frac{1}{2}}}{\Omega} & \text{if } D_0 \text{ is relatively prime to } 4N', \\ 0, & \text{otherwise.} \end{cases}$$

If $g(z)$ is good, then the conclusion of the corollary follows from the theorem.

Let $p$ be a prime that does not divide $4N'$ and that does not split in any imaginary quadratic field having discriminant dividing $4N'$. By [Theorem B($i$), F-H], there is a fundamental discriminant $D'$ such that $p|D'$, $(-1)^k D' > 0$, $D'$ is relatively prime to $4N'$, and $L(f_{D'}, k) \neq 0$. It follows from (11) that $c(D') \neq 0$. For sufficiently large integers $\nu$, the $D'$th Fourier coefficient of the cusp form $G_\nu(z) = g(z) \cdot \theta_{N'^\nu,0,1,1}(z) \in S_{k+1}(4N'^{\nu+1})$ is equal to $c(D')$ and therefore is non-zero. Since $p$ does not divide $4N'$, the new part of $G_\nu(z)$ is non-trivial. Furthermore, the choice of $p$ implies that the new part of $G_\nu(z)$ is not contained in $S_{k+1}^{cm}(4N'^{\nu+1})$ (as the $p$th Fourier coefficient of any cusp form in $S_{k+1}^{cm}(4N'^{\nu+1})$ is zero). This proves that $g(z)$ is good and completes the proof of the corollary.

<div align="right">Q.E.D.</div>

*Proof of Corollary 2.* Let $E/\mathbb{Q}$ be a modular elliptic curve, and let $N$ be the conductor of $E$. For each fundamental discriminant $D < 0$, let $E_D$ be the $D$-quadratic twist of $E$. Let $\omega_D$ be a Neron differential on $E_D$, and let

$$\Omega_D = \int_{E_D(\mathbb{R})} |\omega_D|.$$

Write $\Omega_E$ for $\Omega_{-4}$. Then

$$(12) \qquad \Omega_D D_0^{\frac{1}{2}} = \Omega_E.$$

Since $E$ is modular, it follows that

$$(13) \qquad \frac{L(E_D, 1)}{\Omega_D} \in \mathbb{Q}.$$

However, much more is conjectured to be true. The Birch and Swinnerton-Dyer Conjecture states that if $L(E_D, 1) \neq 0$, then

$$(14) \qquad \frac{L(E_D, 1)}{\Omega_D} = \frac{|\text{Ш}(E_D)|}{|E_D(\mathbb{Q})_{tor}|^2} Tam(E_D),$$

where $Tam(E_D)$ is the Tamagawa factor. Moreover, $Tam(E_D)/Tam(E_{-4})$ is an integer divisible only by the primes 2 or 3. By the work of Mazur [Ma] the order of the torsion group $E_D(\mathbb{Q})_{tor}$ is not divisible by any prime $p > 7$. Since $E$ is modular, there is a newform $f(z) \in S_2(N, \chi_0)$ such that $L(E, s) = L(f, s)$. More generally,

$$(15) \qquad L(E_D, s) = L(f_D, s).$$

It follows from (12) and (13) that $\Omega_E$ is a period for $f(z)$. The result is now a straight-forward consequence of Corollary 1.

<div align="right">Q.E.D.</div>

*Proof of Corollary 3.* Rubin [Ru] has shown that if $E/\mathbb{Q}$ is an elliptic curve having complex multiplication by the CM field $K$, then the only possible primes $p \nmid |O_K^\times|$ dividing $|\text{Ш}(E_D)|$ are those predicted by the Birch and Swinnerton-Dyer Conjecture. This, together with (12), (13), and (14) implies that

$$\ell > 7, \quad |\text{Ш}(E_D)| \equiv 0 \pmod{\ell} \quad \Rightarrow \quad \ell \text{ divides the numerator of } \frac{L(E_D, 1)\sqrt{D_0}}{\Omega_E}.$$

<div align="right">Q.E.D.</div>

There has been some interest in obtaining results for $L(f_D, k)$ where the discriminants $D$ satisfy prescribed local conditions. For instance Kolyvagin's theorem [Ko1] for modular elliptic curves $E$ with conductor $N$ required the existence of at least one fundamental discriminant $D < 0$ for which $L(E_D, s)$ has a simple pole at $s = 1$ with the additional property that $\chi_D(p) = 1$ for every prime $p | N$. This condition has been proved in many ways in the works of Bump-Friedberg-Hoffstein, Iwaniec, and Murty-Murty (see [B-F-H, I, M-M]). More recently, these results have been extended by Friedberg and Hoffstein [F-H]. In Corollary 4 we obtain mod $\ell$ non-vanishing results where the fundamental discriminants similarly satisfy prescribed local conditions.

*Proof of Corollary 4.* If $g(z) = \sum_{n=1}^{\infty} c(n)q^n$ is the cusp form given in the proof of Corollary 1, then for any $0 \leq r < t$,

$$g_{r,t}(z) := \sum_{n \equiv r \pmod{t}} c(n)q^n$$

is a cusp form with level $N't^2$. Consequently it is easy to define a cusp form $g^*(z)$ whose non-zero coefficients are supported by fundamental discriminants satisfying the prescribed conditions. The rest of the argument is the same as the one establishing Corollary 1.

<div align="right">Q.E.D.</div>

*Proof of Corollary 5.* The assertion that $E$ has rank 1 is well known and follows from the works of Gross-Zagier and Kolyvagin [Gr-Z, Ko1]. The rest of the result is proved in [Corollary G, Ko2] conditional on part $(i)$ of Kolyvagin's conjecture where the discriminants satisfy the congruence condition $D \equiv \square \pmod{4N}$. Corollary 4 removes this condition for all but finitely many primes $\ell$, and the proof of Theorem 1 shows how to effectively determine this finite set (see Remark 4).

<div align="right">Q.E.D.</div>

**Remark 1.** By a theorem of Kolyvagin [Ko1] the curves $E_D$ in Corollaries 2 and 3 have rank 0. Moreover by the works of Diamond, Kramer, Shimura, Taylor, and Wiles, the set of modular elliptic curves includes curves having complex multiplication, curves which are semistable at the primes 3 and 5, and also curves with full rational $2-$torsion [Di, Di-K, Sh2, Ta-Wi, Wi].

**Remark 2.** Theorem 1 is a result about the coefficients of an eigenform $g(z)$ of half-integral weight. However, the proof can be applied in a slightly different setting. Let $N$ be an odd square-free integer, and let $f(z) \in S_{2k}(N, \chi_0)$ be a newform. Kohnen and Zagier [K], [K-Z] have constructed an explicit cusp form $g(z) = \sum_{n=1}^{\infty} c(n)q^n \in S_{k+\frac{1}{2}}(N)$ for which $c(n) = 0$ unless $(-1)^k n \equiv 0, 1 \pmod 4$ and for which

$$L(f_D, k) = 2^{-\nu(N)} |D|^{\frac{1}{2}-k} \frac{(k-1)!}{\pi^k} \frac{\langle f, f \rangle}{\langle g, g \rangle} |c(|D|)|^2$$

for any fundamental discriminant $D$ for which $(-1)^k D > 0$ and $\chi_D(\ell) = w_\ell$, the eigenvalue of the Atkin-Lehner involution at $\ell$, for each prime $\ell$ dividing $N$. The conclusion and proof of Theorem 1 apply (*mutatis mutandis*).

**Remark 3.** In view of Corollary 4 it is easy to see that Corollaries 2 and 3 have obvious analogs where local conditions are prescribed. For brevity we omit stating them explicitly.

**Remark 4.** Let $E/\mathbb{Q}$ be a modular elliptic curve with conductor $N$ for which $L(E, s)$ has a simple zero at $s = 1$. Let $g(z)$ be the relevant half-integral weight eigenform, and let $h(z)$ be a theta function for which $G_h(z) = g(z) \cdot h(z)$ satisfies $(G1)$ and $(G2)$. Let $S_1$ denote the set of primes $\ell$ not satisfying $(L1) - (L6)$. If $c(E)$ denotes the Manin constant for $E$, and $D$ denotes the discriminant of $End(E)$, then define $S_2$ by

$$S_2 := \{\ell \mid 6D, \ell \mid Tam(E), \ell \mid c(E)\}.$$

If $E$ does not have complex multiplication, then let $S_3$ denote the finite set of primes $\ell$ for which the $\ell-$adic representation of the Tate module is not surjective, and if $E$ has complex multiplication then let $S_3$ be empty. The conclusion of Corollary 5 holds for every prime $\ell \notin S_1 \cup S_2 \cup S_3$.

**Remark 5.** Our method shows that if there is a single suitable quadratic twist $f_{D'}$ for which $L(f_{D'}, k) \neq 0$ (i.e. 'being good'), then for all but finitely many primes $\ell$ there are infinitely many fundamental $D$ for which the 'algebraic part' of $L(f_D, k) \not\equiv 0 \pmod{\ell}$. Jochnowitz's main results are similarly stated [Jo]. To obtain clean statements we used recent results of Friedberg and Hoffstein [F-H] guaranteeing the existence of infinitely many such twists.

## 3. Examples

**Example 1.** If $N \neq 0$ is an integer, then let $E(N)$ denote the elliptic curve/$\mathbb{Q}$

$$E(N): \quad y^2 = x^3 - N^2 x.$$

Let $g_1(z) := \sum_{n=1}^{\infty} a_1(n)q^n \in S_{\frac{3}{2}}(128, \chi_0)$ and $g_2(z) := \sum_{n=1}^{\infty} a_2(n)q^n \in S_{\frac{3}{2}}(128, \chi_2)$ be the eigenforms defined by

$$g_i(z) := \eta(8z)\eta(16z)\Theta(2^i z).$$

Recall that $\eta(z) := q^{1/24} \prod_{n=1}^{\infty}(1 - q^n)$. If $N \geq 1$ is an odd square-free integer, then Tunnell [T] proved, assuming the Birch and Swinnerton-Dyer Conjecture, that

(16) $$|\text{III}(E(iN))| = \left(\frac{a_i(N)}{2^{\nu(N)}}\right)^2 \quad \text{if} \ \ a_i(N) \neq 0.$$

If $F(z) := g_1(z)\Theta(4z) = g_2(z)\Theta(2z) = \eta(8z)\eta(16z)\Theta(2z)\Theta(4z)$, then

$$F^*(z) := \sum_{n \not\equiv 7 \pmod 8} A(n)q^n - \sum_{n \equiv 7 \pmod 8} A(n)q^n$$

is in $S_2^{\text{new}}(128, \chi_0)$ and does not have complex multiplication. The new part of $F(z)$ is easily seen to be a linear combination of $F^*(z)$ and its twists by Dirichlet characters modulo 8, so $g_1(z)$ and $g_2(z)$ are both good. By Theorem 1, for all sufficiently large primes $\ell$, there are infinitely many odd, square-free integers $N$ and $M$ for which $a_1(N) \not\equiv 0 \pmod{\ell}$ and $a_2(M) \not\equiv 0 \pmod{\ell}$. In fact, a quick inspection of the proof of the theorem shows that in this case 'sufficiently large' means that the image of $\rho_{F^*, \ell}$ is conjugate to $GL_2(\mathbb{F}_\ell)$. The eigenform $F^*(z)$ is the newform associated to the elliptic curve $y^2 = x^3 + 2x^2 - x$ which is a quadratic twist of $y^2 = x^3 - 2x^2 - x$. Serre [5.9.1, Se] has shown that the image of the mod $\ell$ Galois representation of the latter curve is $GL_2(\mathbb{F}_\ell)$ for every odd prime $\ell$. The same is therefore true for $\rho_{F^*, \ell}$.

By Rubin's theorem [Ru] and (16), if $\ell$ is an odd prime, then there are infinitely many odd square-free integers $N$ and $M$ for which $E(N)$ and $E(2M)$ have rank 0 and $\text{III}(E(N))$ and $\mathrm{III}(E(2M))$ have no elements of order $\ell$. The analogous statement when $\ell = 2$ is well known and follows from $2$−descents.

**Example 2.** Let $\Delta(z) = \eta^{24}(z) = \sum_{n=1}^{\infty} \tau(n)q^n \in S_{12}(1)$ denote Ramanujan's cusp form, and let $g(z) = \sum_{n=1}^{\infty} c(n)q^n \in S_{\frac{13}{2}}(4, \chi_0)$ denote the eigenform defined by

$$g(z) := \frac{\Theta^9(z)\eta^8(4z)}{\eta^4(2z)} - \frac{18\Theta^5(z)\eta^{16}(4z)}{\eta^8(2z)} + \frac{32\Theta(z)\eta^{24}(4z)}{\eta^{12}(2z)}.$$

Kohnen and Zagier proved [K-Z] that if $D > 0$ is a fundamental discriminant, then

$$(17) \qquad L(\Delta_D, 6) = \left(\frac{\pi}{D}\right)^6 \frac{\sqrt{D}}{5!} \frac{\langle \Delta(z), \Delta(z) \rangle}{\langle g(z), g(z) \rangle} \cdot c(D)^2.$$

If $F(z) := g(z)\Theta(z) \in S_7(4, \chi_{-1})$, then $F(z) = \dfrac{\mathfrak{B}_1(z) + \mathfrak{B}_2(z)}{2}$, where $\mathfrak{B}_1(z)$ and $\mathfrak{B}_2(z)$ are complex conjugate newforms in $S_7(4, \chi_{-1})$ and where $\mathfrak{B}_1(z)$ is given by

$$\mathfrak{B}_1(z) = \left(1 - \frac{\sqrt{-15}}{15}\right)F(z) + \frac{\sqrt{-15}}{30}F(z)|U_2.$$

The first few terms of the Fourier expansion of $\mathfrak{B}_1(z) = \sum_{n=1}^{\infty} b(n)q^n$ are:

$$\mathfrak{B}_1(z) = q + (2 - 2\sqrt{-15})q^2 + 8\sqrt{-15}q^3 - (56 + 8\sqrt{-15})q^4 + 10q^5 + \dots$$

Since $\mathfrak{B}_1(z)$ and $\mathfrak{B}_2(z)$ do not have complex multiplication, for all sufficiently large primes $\ell$ there exist infinitely many fundamental discriminants $D > 0$ for which $c(D) \not\equiv 0 \pmod{\ell}$ (see Remark 2). As $\mathfrak{B}_2(z)$ is the newform associated to the twist of $\mathfrak{B}_2(z)$ by $\chi_{-1}$, an inspection of the proof of Theorem 1 shows that in this case 'sufficiently large' means $\rho_{\mathfrak{B}_1, \ell}$ is 'big' (i.e. $(L6)$ holds). It is straightforward to check that $\rho_{\mathfrak{B}_1, \ell}$ is irreducible if $\ell \neq 2$, 5, or 61. This can be done by writing down all the possibilities for a reducible $\rho_{\mathfrak{B}_1, \ell}$ and then comparing with the Fourier coefficients of $\mathfrak{B}_1(z)$ using $(R2)$. As in the proofs of [Lemma 5.7, R1] and [Theorem 2.1, R3], if $\rho_{\mathfrak{B}_1, \ell}$ is not 'big', then either its image is dihedral, or for each prime $p \nmid N\ell$

$$\overline{b}(p)^2/\chi_{-1}(p)p^6 \in \{0, 1, 2, 4\} \quad \text{or} \quad \overline{b}(p)^4 - 3\chi_{-1}(p)\overline{b}(p)^2p^6 + p^{12} = 0,$$

where $\overline{b}(p)$ denotes the image of $b(p)$ in $\mathbb{F}_{\mathfrak{B}_1, \ell}$. A simple check of the Fourier coefficients of $\mathfrak{B}_1$ shows that $\rho_{\mathfrak{B}_1, \ell}$ is 'big' if $\ell \neq 2, 3, 5$ or 61. We leave it to the reader to determine what happens at $\ell = 2$ and 3. For $\ell = 5$ and $\ell = 61$ we find that

$$\mathfrak{B}_1(z) + \mathfrak{B}_2(z) \equiv \begin{cases} E_{\chi_{-1}\omega, \omega^{-1}} + E_{\omega, \chi_{-1}\omega^{-1}} \pmod{5}, \\ E_{\chi_0, \chi_{-1}} + E_{\chi_{-1}, \chi_0} \pmod{61}, \end{cases}$$

where $E_{\phi, \psi}$ is the Eisenstein series whose $L$-function is $L(\phi, s)L(\psi, s - 6)$. These congruences, together with the definition of $F(z)$, yield the following Kronecker-style congruences

$$(18) \qquad \sum_{k=-\infty}^{\infty} c(N - k^2) \equiv \begin{cases} \frac{1}{2}N \sum_{d|N} \left(\chi_{-1}(d) + \chi_{-1}(N/d)\right)d^4 \pmod{5}, \\ \frac{1}{2} \sum_{d|N} \left(\chi_{-1}(d) + \chi_{-1}(N/d)\right)d^6 \pmod{61}. \end{cases}$$

Generally a form $g(z) = \sum_{n=0}^{\infty} c(n)q^n \in M_{k+\frac{1}{2}}(N)$ with algebraic integer coefficients satisfies a Kronecker congruence mod $\ell$ if there exists some $h(z) \in \Theta(N')$ for which

$$(i) \quad G_h^{new}(z) = \sum_{f \in \mathcal{N}} \alpha_f f(z) \neq 0,$$

$$(ii) \quad \rho_{\ell,f} \text{ is reducible if } \alpha_f \neq 0,$$

$$(iii) \quad |\alpha_f|_\ell \leq 1 \quad \text{if } \alpha_f \neq 0.$$

**Example 3.** In this example we examine $\ell$-indivisibility of class numbers of imaginary quadratic fields. Since obtaining indivisibility results follow easily by Kronecker's class number relation [Ha], we select an arithmetic progression of discriminants for which these indivisibility results do not follow so easily. We investigate the $\ell$-indivisibility of the class numbers $h(-32n - 20)$. Using the following identity of Jacobi

$$\frac{\eta^2(16z)}{\eta(8z)} := \sum_{n=0}^{\infty} q^{(2n+1)^2},$$

it is easy to see that

$$\frac{\eta^4(32z)}{\eta(8z)} = \frac{\eta^2(16z)}{\eta(8z)} \left( \frac{\eta^2(32z)}{\eta(16z)} \right)^2 = \sum_{n \geq 5 \text{ odd}} c(n)q^n = \sum_{x,y,z \geq 1 \text{ odd}} q^{x^2 + 2y^2 + 2z^2}.$$

Since the ternary quadratic form $x^2 + 2y^2 + 2z^2$ is the only class in its genus, if $8n + 5$ is square-free, then $h(-32n - 20) = 2c(n)$. Incidentally, if $n \equiv 5 \pmod 8$, then $c(n)$ is the number of $4-$core partitions of $\frac{n-5}{8}$ (see [O-S]).

If $F(z) \in S_3(32, \chi_{-1})$ is defined by

$$F(z) := \frac{\eta^4(16z)}{\eta(4z)} \eta^3(4z) = q^3 - 2q^7 - q^{11} + 2q^{15} + \dots,$$

then by Jacobi's triple product identity

$$\eta^3(8z) := \sum_{n=0}^{\infty} (-1)^n (2n+1) q^{(2n+1)^2},$$

we find that

$$F(z) = \sum_{\substack{n \geq 5 \\ x \geq 0}} (-1)^x (2x+1) c(n) q^{(n+(2x+1)^2)/2}.$$

Moreover, $F(z) = \dfrac{\mathfrak{B}_2(z) - \mathfrak{B}_1(z)}{8i}$, where $\mathfrak{B}_1(z)$ and $\mathfrak{B}_2(z) \in S_3(32, \chi_{-1})$ are complex conjugate newforms in $S_3(32, \chi_{-1})$ which do not have complex multiplication, and $\mathfrak{B}_1(z)$ is given by

$$\mathfrak{B}_1(z) := F(z)|T_3 - 4iF(z) = q - 4iq^3 + 2q^5 + 8iq^7 - 7q^9 + \cdots$$

Just as in Example 2, $\mathfrak{B}_2(z)$ is also the newform associated to the twist of $\mathfrak{B}_1(z)$ by $\chi_{-1}$. While $\dfrac{\eta^4(32z)}{\eta(8z)}$ is an eigenform, $\dfrac{\eta^4(16z)}{\eta(4z)}$ is not, so we cannot appeal directly to Theorem 1. However the methods used to prove Theorem 1 show that if $\ell$ is any odd prime for which $\rho_{\mathfrak{B}_1,\ell}$ is 'big' (i.e. for which $(L6)$ holds), then there are infinitely many fundamental discriminants $-32n - 20$ for which $h(-32n - 20) \not\equiv 0 \pmod \ell$. The arguments employed in the previous example show that if $\ell \geq 5$, then $\rho_{\mathfrak{B}_1,\ell}$ is big. Again, we leave the case where $\ell = 3$ to the reader.

## 4. Concluding remarks

It is somewhat annoying that Theorem 1 only pertains to good forms $g(z)$. This condition is expected to be very mild. For example, suppose $g(z) = \sum_{n=1}^{\infty} c(n)q^n \in S_{k+\frac{1}{2}}(N, \chi_d)$ is a 'bad' eigenform lifting to a newform $f(z) \in S_{2k}(N', \chi_0)$. By [Cor. 2, Wal] there is at least one arithmetic progression $r \pmod{t}$ for which every square-free positive integer $n \equiv r \pmod{t}$ has the property that the sign of the functional equation of $L(f_{(-1)^k dn}, s)$ is $+1$, and

$$L(f_{(-1)^k dn}, k) = c(n)^2 \cdot A_n$$

where $A_n$ is an explicit non-zero constant. By hypothesis, for every $\nu \geq 1$ $g(z)\Theta(N^\nu z) \in S_{2k}^{cm}(4N^\nu)$. Therefore it turns out that $c(n) = 0$ for every 'inert' $n$, a set of positive integers with arithmetic density 1, so $L(f_{(-1)^k dn}, k) = 0$ for almost every square-free $n \equiv r \pmod{t}$. Since it is widely believed that there is no such newform $f(z)$, we are lead to the following conjecture.

**Conjecture.** If $g(z) \in M_{k+\frac{1}{2}}(N) \setminus \Theta(N)$, then $g(z)$ is good.

## Acknowledgements

The authors thank B. Duke, V. Kolyvagin, D. Lieman, and K. Soundararajan for helpful comments.

School of Mathematics, Institute for Advanced Study, Princeton, New Jersey 08540
    *E-mail address*: ono@math.ias.edu

Department of Mathematics, The Pennsylvania State University, University Park, Pennsylvania 16802
    *E-mail address*: ono@math.psu.edu

Department of Mathematics, Princeton University, Princeton, New Jersey 08540
    *E-mail address*: cmcls@math.princeton.edu